# ON THE VALUE OF OPTIMAL STOPPING GAMES

By Erik Ekström and Stephane Villeneuve[1]

*University of Manchester and University of Toulouse*

We show, under weaker assumptions than in the previous literature, that a perpetual optimal stopping game always has a value. We also show that there exists an optimal stopping time for the seller, but not necessarily for the buyer. Moreover, conditions are provided under which the existence of an optimal stopping time for the buyer is guaranteed. The results are illustrated explicitly in two examples.

**1. Introduction.** In this paper we study a perpetual optimal stopping game between two players, the "buyer" and the "seller." Both players choose a stopping time each, say $\tau$ and $\gamma$, and at the time $\tau \wedge \gamma := \min\{\tau, \gamma\}$, the seller pays the amount

$$Y_1(\tau)\mathbb{1}_{\{\tau \leq \gamma\}} + Y_2(\gamma)\mathbb{1}_{\{\tau > \gamma\}} \tag{1.1}$$

to the buyer. Here $Y_1$ and $Y_2$ are two stochastic processes satisfying $0 \leq Y_1(t) \leq Y_2(t)$ for all $t$ almost surely. Clearly, the seller wants to minimize the amount in (1.1) and the buyer wants to maximize this amount.

We consider discounted optimal stopping games defined in terms of two continuous contract functions $g_1$ and $g_2$ satisfying $0 \leq g_1 \leq g_2$ and a one-dimensional diffusion process $X(t)$. More precisely, given a constant discounting rate $\beta > 0$, let

$$Y_1(t) = e^{-\beta t}g_1(X(t))$$

and

$$Y_2(t) = e^{-\beta t}g_2(X(t)).$$

Define the mapping $R_x$ from the set of pairs $(\tau, \gamma)$ of stopping times to the set $[0, \infty]$ by

$$R_x(\tau, \gamma) := \mathbb{E}_x e^{-\beta \tau \wedge \gamma}(g_1(X(\tau))\mathbb{1}_{\{\tau \leq \gamma\}} + g_2(X(\gamma))\mathbb{1}_{\{\tau > \gamma\}}). \tag{1.2}$$

Received July 2005; revised January 2006.
[1]Supported by the Fonds National de la Science.

Thus, $R_x(\tau,\gamma)$ is the expected discounted pay-off when the players use the stopping times $\tau$ and $\gamma$ as stopping strategies. Here the index $x$ indicates that the diffusion $X$ is started at $x$ at time 0. In (1.2), and in similar situations below, we use the convention that

$$f(X(\sigma)) = 0 \quad \text{on } \{\sigma = \infty\},$$

where $f$ is a function and $\sigma$ is a random time. Next define the lower value $\underline{V}$ and the upper value $\overline{V}$ as

$$\underline{V}(x) := \sup_{\tau} \inf_{\gamma} R_x(\tau,\gamma)$$

and

$$\overline{V}(x) := \inf_{\gamma} \sup_{\tau} R_x(\tau,\gamma),$$

respectively, where the supremums and the infimums are taken over random times $\tau$ and $\gamma$ that are stopping times. It is clear that

$$g_1(x) \leq \underline{V}(x) \leq \overline{V}(x) \leq g_2(x)$$

(the first and the last inequality follow from choosing $\tau = 0$ or $\gamma = 0$ in the definitions of $\underline{V}$ and $\overline{V}$, resp.). If, in addition, the inequality

$$\underline{V}(x) \geq \overline{V}(x)$$

holds, that is, if $\underline{V}(x) = \overline{V}(x)$, then the stochastic game is said to have a value. In such cases, we denote the common value $\underline{V}(x) = \overline{V}(x)$ by $V(x)$. If there exist two stopping times $\tau'$ and $\gamma'$ such that

(1.3) $$R_x(\tau,\gamma') \leq R_x(\tau',\gamma') \leq R_x(\tau',\gamma)$$

for all stopping times $\tau$ and $\gamma$, then the pair $(\tau',\gamma')$ is referred to as a saddle point for the stochastic game. It is clear that if there exists a saddle point for the stochastic game, then the game also has a value.

It is well known (compare [2, 3, 10, 11, 13] and [15]) that under the integrability condition

(1.4) $$\mathbb{E}_x \left( \sup_{0 \leq t < \infty} e^{-\beta t} g_2(X(t)) \right) < \infty$$

and the condition

$$\lim_{t \to \infty} e^{-\beta t} g_2(X(t)) = 0,$$

the stochastic game has a value $V$. Moreover, the two stopping times

(1.5) $$\tau^* := \inf\{t : V(X(t)) = g_1(X(t))\}$$

and

(1.6) $$\gamma^* := \inf\{t : V(X(t)) = g_2(X(t))\}$$

together form a saddle point for the game. Below we prove the existence of a value under no integrability conditions at all. To do this, we use the connection between excessive functions and concave functions; compare [6] and [7]. More specifically, using concave functions, we produce a candidate $V^*$ for the value function, and then we prove that $\underline{V} \geq V^* \geq \overline{V}$. Thus, there exists a value of the game, and this value is given by the candidate function $V^*$. One should note that we prove the existence of a value for perpetual optimal stopping games, that is, when there is no upper bound on the stopping times $\tau$ and $\gamma$. It remains an open question if all optimal stopping games with a finite time horizon have values.

One easily finds examples of stochastic differential games where the pair $(\tau^*, \gamma^*)$ of stopping times defined by (1.5) and (1.6) is not a saddle point; compare, for instance, the examples in Section 5.1. We prove below, however, that $\gamma^*$ is always optimal for the seller. More precisely, we deal with the following concepts closely related to the notion of a saddle point: a stopping time $\tau'$ is optimal for the buyer if

$$R_x(\tau', \gamma) \geq \overline{V}(x)$$

for all stopping times $\gamma$, and a stopping time $\gamma'$ is optimal for the seller if

$$R_x(\tau, \gamma') \leq \underline{V}(x)$$

for all stopping times $\tau$. Note that

$\tau'$ is optimal for the buyer and $\gamma'$ is optimal for the seller

$\iff \quad (\tau', \gamma')$ is a saddle point.

Also note that if $\tau'$ is optimal for the buyer, then

$$\overline{V}(x) \leq \inf_\gamma R_x(\tau', \gamma) \leq \underline{V}(x) \leq \overline{V}(x),$$

so the game has a value $V(x)$ which is given by

$$V(x) = \inf_\gamma R_x(\tau', \gamma).$$

Similarly, if $\gamma'$ is optimal for the seller, then the existence of a value $V(x)$ follows, and

$$V(x) = \sup_\tau R_x(\tau, \gamma').$$

The outline of the paper is as follows. In Section 2 we specify the assumptions on the diffusion $X$ and we show that a stochastic game with an infinite time horizon always has a value. This is done without the integrability condition (1.4); compare Theorem 2.5. We also show that $\gamma^*$ is an optimal stopping time for the seller. The method used in the proof of Theorem 2.5 also gives a characterization of the value function in terms of concave functions. As a straightforward consequence of this characterization, the smooth-fit principle is deduced in Section 3. In Section 4 we provide additional conditions under which $\tau^*$ is optimal for the buyer, that is, $(\tau^*, \gamma^*)$ is a saddle point. Finally, in Section 5 we explicitly determine the value of two different game options, both of which may be regarded game versions of the American call option. In these examples, the integrability condition (1.4) is not fulfilled, so they are not covered by the theory in previous literature.

**2. The value of a stochastic differential game.** Let $X$ be a stochastic process with dynamics

$$(2.1) \qquad dX(t) = \mu(X(t))\,dt + \sigma(X(t))\,dW(t),$$

where $\mu$ and $\sigma$ are given functions and $W$ is a standard Brownian motion. We assume that the two end-points of the state space of $X$ are 0 and $\infty$, and we assume for simplicity that both these end-points are natural. We also assume that the functions $\mu(\cdot)$ and $\sigma(\cdot)$ are continuous and that $\sigma(x) > 0$ for all $x \in (0, \infty)$. It follows that the equation (2.1) has a (weak) solution which is unique in the sense of probability law; see Chapter 5.5 in [12]. Moreover, $X$ is a regular diffusion, that is, for all $x, y \in (0, \infty)$, we have that $y$ is reached in finite time with a positive probability if the diffusion is started from $x$.

The second-order ordinary differential equation

$$(2.2) \qquad \mathcal{L}u(x) := \frac{\sigma^2(x)}{2} u_{xx} + \mu(x) u_x - \beta u = 0$$

has two linearly independent solutions $\psi, \varphi : (0, \infty) \to \mathbb{R}$ which are uniquely determined (up to multiplication with positive constants) by requiring one of them to be positive and strictly increasing and the other one to be positive and strictly decreasing; compare [5]. We let $\psi$ be the increasing solution and $\varphi$ the decreasing solution. Since 0 and $\infty$ are assumed to be natural boundaries of $X$, we have $\psi(0+) = 0 = \varphi(\infty)$. We also let $F : (0, \infty) \to (0, \infty)$ be the strictly increasing positive function defined by

$$F(x) := \frac{\psi(x)}{\varphi(x)}.$$

Recall that a function $u : (0, \infty) \to \mathbb{R}$ is said to be $F$-concave in an interval $J \subset (0, \infty)$ if

$$u(x) \geq u(l) \frac{F(r) - F(x)}{F(r) - F(l)} + u(r) \frac{F(x) - F(l)}{F(r) - F(l)}$$

for all $l, x, r \in J$ with $l < x < r$. Equivalently, the function $u(F^{-1}(\cdot))$ is concave. $F$-convexity of a function is defined similarly.

Below we use the following two theorems relating concave and convex functions to the value functions of optimal stopping problems. The first one is Proposition 4.2 in [6]. The proof of the second one follows along the lines of the proofs of Propositions 3.2 and 4.2 in [6] and is therefore omitted.

THEOREM 2.1. *Let $l, r$ be such that $0 < l < r < \infty$, let $g : [l, r] \to [0, \infty)$ be measurable and bounded, and let*

$$U(x) := \sup_{\tau \leq \tau_{l,r}} \mathbb{E}_x e^{-\beta \tau} g(X(\tau)),$$

*where*

$$\tau_{l,r} := \inf\{t : X(t) \notin (l, r)\}.$$

*Then $U$ is the smallest majorant of $g$ such that $U/\varphi$ is $F$-concave on $[l, r]$.*

THEOREM 2.2. *Let $l, r$ be such that $0 < l < r < \infty$, let $g : [l, r] \to [0, \infty)$ be measurable and bounded, and let*

$$U(x) := \inf_{\gamma \leq \gamma_{l,r}} \mathbb{E}_x e^{-\beta \gamma} g(X(\gamma)),$$

*where*

$$\gamma_{l,r} := \inf\{t : X(t) \notin (l, r)\}.$$

*Then $U$ is the largest minorant of $g$ such that $U/\varphi$ is $F$-convex on $[l, r]$.*

REMARK. Note that it is important in Theorem 2.2 that the stopping times $\gamma$ are to be chosen among stopping times not exceeding the first exit time $\gamma_{l,r}$ of $X(t)$ from the interval $(l, r)$. If, for example, the choice $\gamma = \infty$ would be included, then $U$ would be identically 0.

Below we find our candidate value function $V^*$ in the set

$$\mathbb{F} = \{f : (0, \infty) \to [0, \infty) : f \text{ is continuous, } g_1 \leq f \leq g_2,$$
$$f/\varphi \text{ is } F\text{-concave in every interval in which } f < g_2\}.$$

Note that $\mathbb{F}$ is nonempty since $g_2 \in \mathbb{F}$. We work below with the functions $H_i : (0, \infty) \to [0, \infty)$, $i = 1, 2$, defined by

(2.3) $$H_i(y) := \frac{g_i(F^{-1}(y))}{\varphi(F^{-1}(y))}$$

and the set

$$\mathbb{H} = \{h : (0, \infty) \to [0, \infty) : h \text{ is continuous, } H_1 \leq h \leq H_2,$$
$$h \text{ is concave in every interval in which } h < H_2\}.$$

Note that the functions in $\mathbb{F}$ are precisely the functions $\varphi \cdot (h \circ F)$ for some function $h \in \mathbb{H}$.

LEMMA 2.3. *Let $\{h_n\}_{n=1}^\infty$ be a sequence of functions in $\mathbb{H}$. Then the function $h$ defined by*

$$h(y) := \inf_n h_n(y)$$

*is an element of $\mathbb{H}$.*

PROOF. First we claim that the minimum of two functions in $\mathbb{H}$ is again in $\mathbb{H}$. To see this, assume that $h_1, h_2 \in \mathbb{H}$ and let $h := h_1 \wedge h_2$. Clearly, $h$ is continuous and satisfies $H_1 \leq h \leq H_2$. Let $y \in (0, \infty)$ satisfy $h(y) < H_2(y)$. Consider the two separate cases $h_1(y) \neq h_2(y)$ and $h_1(y) = h_2(y) < H(y)$. In the first case, there exists an open interval containing $y$ such that $h = h_1$ or $h = h_2$ in this interval and, thus, $h$ is concave in this interval. For the second case, there exists an open interval containing $y$ such that both $h_1$ and $h_2$ are concave. Since the minimum of two concave functions is concave, $h$ is also concave in this interval. It follows that $h$ is concave in every interval in which $h < H_2$, which shows that $h \in \mathbb{H}$.

Thus, we may, without loss of generality, assume that $h_{n+1} \leq h_n$ for all $n$. Let $h(y) := \inf_n h_n(y)$ and define

$$U := \{y : h(y) < H_2(y)\}.$$

Note that $h$, being the infimum of continuous functions, is upper semi-continuous, so $U$ is open. Choose two points $l, r \in U$ with $l < r$ and $[l, r] \subset U$. The interval $[l, r]$ is compact, and it is covered by the increasing family $\{U_n\}_{n=1}^\infty$ of open sets

$$U_n := \{y : h_n(y) < H_2(y)\}.$$

Hence, there exists an integer $N$ such that $[l, r] \subset U_n$ for all $n \geq N$. For such $n$, $h_n$ is concave on $[l, r]$, and therefore, also $h$ is concave on this interval. Consequently, $h$ is concave on each interval contained in $U$, and thus also continuous at all points in $U$.

To show that $h \in \mathbb{H}$, it remains to check that $h$ is continuous also at all boundary points of $U$. Let $l \in \overline{U} \setminus U$, where $\overline{U}$ is the closure of $U$ in $(0, \infty)$, and let $\{l_k\}_{k=1}^\infty$ be a sequence of points in $U$ converging to $l$ from the right (left-continuity is dealt with similarly). Because $h$ is upper semi-continuous, it is enough to prove that $h(l) \leq h(l+)$.

Assume first that $(l, l + \varepsilon_0) \subset U$ for some $\varepsilon_0 > 0$. We assume, to reach a contradiction, that there exists $\varepsilon > 0$ such that $h(l) - \varepsilon > h(l+)$. Then there exists $\delta > 0$ such that the straight line $L$ connecting the points $(l, h(l) - \varepsilon)$ and $(l+\delta, h(l+\delta))$ satisfies $h(y) < L(y) < H_2(y)$ for $y \in (l, l+\delta)$. Now, choose a $y \in (l, l + \delta)$. Then there exists an $n$ such that $h_n(y) < L(y)$. For this $n$, $h_n(l) \leq L(l)$ since $h_n$ is concave and $h_n(l + \delta) \geq L(l + \delta)$. Consequently, $h(l) \leq h_n(l) \leq L(l) = h(l) - \varepsilon$, which is the required contradiction.

On the other hand, if there does not exist an $\varepsilon_0 > 0$ such that $h < H_2$ in $(l, l + \varepsilon_0)$, then the previous case can be applied to deduce right-continuity of $h$ at $l$. Indeed, for $\varepsilon > 0$, choose $\delta > 0$ such that

$$|H_2(y) - H_2(l)| \leq \varepsilon \quad \text{for } y \in [l, l + \delta].$$

Without loss of generality, we may assume that $H_2 = h$ at $l + \delta$ (since points with $H_2 = h$ exist arbitrarily close to $l$). Now, for a point $l_k \in (l, l + \delta)$, there exists a maximal (possibly empty) surrounding interval in which $h < H_2$. We know from above that $h$ is concave in the closure of this interval, and thus, $h \geq H_2(l) - \varepsilon$ in the interval. In particular, $h(l_k) \geq H_2(l) - \varepsilon$. Since we also have $h(l_k) \leq H_2(l_k) \leq H_2(l) + \varepsilon$, and since $\varepsilon$ is arbitrary, it follows that $h(l_k) \to H_2(l) = h(l)$ as $k \to \infty$. Hence, $h$ is continuous at $l$, and thus, we have shown that $h \in \mathbb{H}$. $\square$

LEMMA 2.4. *There exists a smallest element $V^* \in \mathbb{F}$. Moreover, the function $V^*/\varphi$ is $F$-convex in every interval in which $V^* > g_1$.*

PROOF. Since the functions in $\mathbb{F}$ are precisely the functions $\varphi(x)h(F(x))$ for some function $h \in \mathbb{H}$, it suffices to show that there exists a smallest element in $\mathbb{H}$ and that this smallest element is convex in every interval of strict majorization of $H_1$. In order to do this, define

$$W(y) := \inf_{h \in \mathbb{H}} h(y).$$

Being the infimum of continuous functions, $W$ is itself upper semi-continuous. Let $\{y_k\}_{k=1}^{\infty}$ be a dense sequence of points in $(0, \infty)$, and for each $k$, let

$$\{h_n^k\}_{n=1}^{\infty} \subseteq \mathbb{H}$$

be a sequence of functions in $\mathbb{H}$ such that $\inf_n h_n^k(y_k) = W(y_k)$. Next, define the function $W^*$ by

$$W^*(y) = \inf_k \inf_n h_n^k(y).$$

According to Lemma 2.3, $W^* \in \mathbb{H}$. Moreover, the nonnegative function $W^* - W$ is lower semi-continuous and vanishes on a dense subset of $(0, \infty)$. It follows that $W \equiv W^*$, so $W \in \mathbb{H}$, which finishes the first part of the proof.

To show the convexity on each interval in which $W > H_1$, let $I$ be such an interval and fix $y' \in I$. By continuity of $H_1$, $H_2$ and $W$, we can find $\delta > 0$ so that

$$\inf_{y \in I^\delta} W(y) \geq \sup_{y \in I^\delta} H_1(y),$$

where $I^\delta := [y' - \delta, y' + \delta]$. Now assume, to reach a contradiction, that there exist points $y_1, y_2 \in I^\delta$ with $y_1 < y' < y_2$ and

$$(2.4) \quad W(y') > W(y_1)\frac{y_2 - y'}{y_2 - y_1} + W(y_2)\frac{y' - y_1}{y_2 - y_1} =: L(y').$$

Since $W$ is continuous, $W(y) > L(y)$ for $y$ in an open set containing $y'$. Let us introduce

$$y_1' = \sup\{y \in [y_1, y'], W(y) = L(y)\}$$

and

$$y_2' = \inf\{y \in [y', y_2], W(y) = L(y)\}.$$

It is now straightforward to check that the function

$$h(y) := \begin{cases} L(y), & \text{if } y \in [y_1', y_2'], \\ W(y), & \text{if } y \notin (y_1', y_2'), \end{cases}$$

satisfies $h \in \mathbb{H}$. However, $h < W$ in $y \in (y_1', y_2')$ contradicts the minimality of $W$, and thus, (2.4) is not true. This means that $W$ is convex at the point $y'$, so, by continuity, $W$ is convex on $I$, which finishes the second part of the proof. □

THEOREM 2.5. *For any starting point $x > 0$, the perpetual optimal stopping game has a value $V(x) := \underline{V}(x) = \overline{V}(x)$. Moreover, $V \equiv V^*$, where $V^*$ is the function appearing in Lemma 2.4, and the stopping time*

$$\gamma^* := \inf\{t : V(X(t)) = g_2(X(t))\}$$

*is an optimal stopping time for the seller.*

PROOF. Let $V^*$ be the function in Lemma 2.4, and choose $x \in (0, \infty)$. To prove the existence of a value, we will show that

$$(2.5) \quad \overline{V}(x) \leq V^*(x) \leq \underline{V}(x).$$

To prove the first inequality, assume that the maximal interval containing $x$ in which $V^* < g_2$ is $(l, r)$ for some points $l < r$ [if $V^*(x) = g_2(x)$, then the first inequality obviously holds since $\overline{V} \leq g_2$]. Assume also, for the moment, that $0 < l$ and $r < \infty$. It follows that $V^*(l) = g_2(l)$ and $V^*(r) = g_2(r)$. Inserting $\gamma = \gamma_{l,r}$ in the definition of $\overline{V}$ yields

$$\overline{V}(x) \leq \sup_\tau \mathbb{E}_x e^{-\beta \tau \wedge \gamma_{l,r}}(g_1(X(\tau))\mathbb{1}_{\{\tau \leq \gamma_{l,r}\}} + g_2(X(\gamma_{l,r}))\mathbb{1}_{\{\tau > \gamma_{l,r}\}})$$

$$(2.6) \quad \leq \sup_{\tau \leq \gamma_{l,r}} \mathbb{E}_x e^{-\beta \tau} g^*(X(\tau))$$

$$=: U^*(x),$$

where the function $g^*$ is defined by

$$(2.7) \qquad g^*(x) = \begin{cases} g_1(x), & \text{if } x \in (l,r), \\ g_2(x), & \text{if } x \in \{l,r\}. \end{cases}$$

Note that $V^*$ majorizes $g^*$ and that $V^*/\varphi$ is $F$-concave on $(l,r)$. According to Theorem 2.1, $U^*$ is the smallest such function, so $U^*(x) \leq V^*(x)$. Consequently,

$$(2.8) \qquad \overline{V}(x) \leq V^*(x).$$

Now, if we instead have $0 = l$ and/or $r = \infty$, then the above reasoning again applies if we plug in $\gamma_r := \inf\{t : X(t) \geq r\}$, $\gamma_l := \inf\{t : X(t) \leq l\}$ or $\gamma = \infty$ in the definition of $\overline{V}$ and use Propositions 5.3 or 5.11 in [6] instead of Theorem 2.1.

To show the second inequality in (2.5), we argue similarly. Choose an $x$ and let $(l,r)$ be a maximal interval containing $x$ in which $V^* > g_1$. As above, let us first assume that

$$(2.9) \qquad 0 < l \quad \text{and} \quad r < \infty.$$

Inserting $\tau = \tau_{l,r}$ in the definition of $\underline{V}$ gives

$$\underline{V}(x) \geq \inf_\gamma \mathbb{E}_x e^{-\beta \tau_{l,r} \wedge \gamma}(g_1(X(\tau_{l,r}))\mathbb{1}_{\{\tau_{l,r} \leq \gamma\}} + g_2(X(\gamma))\mathbb{1}_{\{\tau_{l,r} > \gamma\}})$$

$$= \inf_{\gamma \leq \tau_{l,r}} \mathbb{E}_x e^{-\beta \gamma} g_*(X(\gamma)),$$

where the function $g_*$ is given by

$$g_*(x) = \begin{cases} g_2(x), & \text{if } x \in (l,r), \\ g_1(x), & \text{if } x \in \{l,r\}. \end{cases}$$

Thus, since $V^*/\varphi$ is $F$-convex in $(l,r)$ (see Lemma 2.4), it follows from Theorem 2.2 that $\underline{V}(x) \geq V^*(x)$. Thus, we have shown the second inequality in (2.5) under the assumption (2.9).

Now, if (2.9) is not the case, then the second inequality in (2.5) requires some slightly more involved analysis. For example, assume that

$$(2.10) \qquad 0 < l \quad \text{and} \quad r = \infty.$$

To prove $\underline{V}(x) \geq V^*(x)$, in this case we do not plug in $\tau_l := \inf\{t : X(t) \leq l\}$ in the definition of $\underline{V}$, but we rather use the stopping times $\tau_{l,N} = \inf\{t : X(t) \notin (l,N)\}$ for different $N \geq l$ (compare the remark following the current proof). Thus, for any $N \geq x$, choosing $\tau = \tau_{l,N}$ in the definition of $\underline{V}$ gives

$$(2.11) \quad \underline{V}(x) \geq \inf_\gamma R_x(\tau_{l,N}, \gamma) = \inf_{\gamma \leq \tau_{l,N}} \mathbb{E}_x e^{-\beta \gamma} g_*(X(\gamma)) =: V_N(x),$$

where

$$g_*(x) = \begin{cases} g_2(x), & \text{if } x \in (l,N), \\ g_1(x), & \text{if } x = \{l,N\}. \end{cases}$$

From Theorem 2.2, it follows that $V_N$ is majorized by $g_*$, that $V_N/\varphi$ is $F$-convex on $[l, N]$, and that $V_N$ is the largest function with these properties. It is clear from (2.8) and (2.11) that

$$\sup_{N \geq x} V_N(x) \leq \underline{V}(x) \leq V^*(x).$$

We show below that we in fact have

(2.12) $$\sup_{N \geq x} V_N(x) = V^*(x).$$

Note that (2.12) implies that

$$\underline{V}(x) = V^*(x)$$

and therefore also the existence of a value. To prove (2.12), we will work in the coordinates $y$ defined by $y = F(x)$.

Let $H_i$, $i = 1, 2$, be defined by $H_i = \frac{g_i}{\varphi} \circ F^{-1}$. Then

$$W_{N'} := \frac{V_N}{\varphi} \circ F^{-1} : [l', N'] \to \mathbb{R}$$

is the largest convex function majorized by the function

(2.13) $$H(y) := \begin{cases} H_2(y), & \text{if } y \in (l', N'), \\ H_1(y), & \text{if } y \in \{l', N'\}, \end{cases}$$

where $l' := F(l)$ and $N' := F(N)$. Let $W := \frac{V^*}{\varphi} \circ F^{-1}$ (thus, $W$ is the function defined in the proof of Lemma 2.4). The conditions $0 < l$ and $r = \infty$ translate to $l' > 0$, $W(l') = H_1(l')$ and $W(y) > H_1(y)$ for all $y > l'$. Next, for $y > l'$, define

$$\hat{W}(y) := \sup_{N' \geq y} W_{N'}(y).$$

We need to show that $\hat{W} \geq W$. To do this, note that since $W > H_1$ in the interval $[l', \infty)$, we know from Lemma 2.4 that $W$ is convex in this interval. Choose $y_0 > l'$, let

$$k := \lim_{\varepsilon \searrow 0} \frac{W(y_0 + \varepsilon) - W(y_0)}{\varepsilon}$$

be the right derivative of $W$ at $y_0$, and let $L(y) = k(y - y_0) + W(y_0)$ be the steepest tangential of $W$ at $y_0$. Note that $L(y) \leq W(y) \leq H_2(y)$. Now we consider two cases.

First, assuming the existence of a point $N' > y_0$ such that $L(N') = H_1(N')$, the function

$$h(y) = \begin{cases} W(y), & \text{if } y \in [l', y_0], \\ L(y), & \text{if } y \in [y_0, N'], \end{cases}$$

is convex and dominated by $H_2$ in $(l', N')$ and by $H_1$ at the points $l'$ and $N'$. Therefore, $h \leq W_{N'}$ by Theorem 2.2, so

$$\hat{W}(y_0) \geq W(y_0).$$

Second, assume that there is no point $N' > y_0$ such that $L(N') = H_1(N')$. Note that the function

$$h(y) = \begin{cases} W(y), & \text{if } y \in (0, y_0], \\ L(y), & \text{if } y \in [y_0, \infty), \end{cases}$$

is an element of the set $\mathbb{H}$. Since $W$ is the smallest function in this set, it follows that we must have $W(y) = L(y)$ for all $y \geq y_0$. Moreover, for each $\varepsilon > 0$, there exists a point of intersection (to the right of $y_0$) between the line $L^\varepsilon(y) := (k - \varepsilon)(y - y_0) + W(y_0)$ and $H_1$ (otherwise a function in $\mathbb{H}$ can be constructed which is strictly smaller than $W$ in some interval). Now, let $z < W(y_0)$, and consider the straight lines through $(y_0, z)$ that are below $W$ in the interval $[l', y_0]$. Let $k'$ be the slope of the largest such straight line (i.e., $k'$ is the smallest possible slope), denote this line by $L'$, and let $y' \in [l', y_0)$ be the largest value for which $W = L'$. Since $W$ is convex in $[l', \infty)$, we have that $k' < k$, and thus, the straight line through $(y_0, W(y_0))$ with slope $k'$ and the function $H_1$ have a point $(N', H_1(N'))$ of intersection for some $N' > y_0$. Let $L''$ be the straight line between the points $(y_0, z)$ and $(N', H_1(N'))$. Then the function which equals $W$ in $[l', y']$, $L'$ in $[y', y_0]$ and $L''$ in $[y_0, N']$ is convex and smaller than the function $H$ defined as in (2.13). Consequently, the corresponding function $W_{N'}$ satisfies $W_{N'}(y_0) \geq z$. Since $z < W(y_0)$ is arbitrary, it follows that $\hat{W}(y_0) \geq W(y_0)$.

Thus, we have shown under the assumption (2.10) that (2.12) holds, implying the second inequality in (2.5). By symmetry, the above argument also applies in the case when $l = 0$ and $r < \infty$. The remaining case, that is, when $l = 0$ and $r = \infty$, can be handled with similar methods (we omit the details).

Finally, since we have shown that the first inequality in (2.6) actually is an equality, it follows that $\gamma^*$ is optimal for the seller. $\square$

REMARK. Note that the function $W$ in the proof of Lemma 2.4 is the smallest function in the set

$$\mathbb{H} = \{h : (0, \infty) \to [0, \infty) : h \text{ is continuous, } H_1 \leq h \leq H_2,$$

$$h \text{ is concave in every interval in which } h < H_2\},$$

whereas, in general, it is not the largest function in the set

$$\{h : (0, \infty) \to [0, \infty) : h \text{ is continuous, } H_1 \leq h \leq H_2,$$

$$h \text{ is convex in every interval in which } h > H_1\}$$

(although $W$ is a member also of this set). This asymmetry of the function $W$ (and the corresponding one for the function $V^*$) may be regarded as the underlying reason for the asymmetry in the proof of the first and the second inequality in (2.5).

REMARK. Let us introduce the perpetual American option value $V_\infty$ associated with the payoff $g_1$, that is,

(2.14) $$V_\infty(x) := \sup_\tau \mathbb{E}_x e^{-\beta\tau} g_1(X(\tau)).$$

Obviously, $V \leq V_\infty$. An immediate consequence of Theorem 2.5 is that the implication

$$V_\infty(x_0) \geq g_2(x_0) \quad \text{for some } x_0 \in (0,\infty) \implies \{x : V(x) = g_2(x)\} \neq \varnothing$$

holds. Indeed, assume that $V_\infty(x_0) \geq g_2(x_0)$ for some $x_0$ and that $V(x) < g_2(x)$ for all $x \in (0,\infty)$. Then $\gamma^* = \infty$, so $V \equiv V_\infty$ by Theorem 2.5. It follows that $V(x_0) \geq g_2(x_0)$, which is a contradiction.

**3. The smooth-fit principle.** In the following proposition, let $H_1$ and $H_2$ be the functions defined in (2.3) and let $W$ be the smallest element in the set $\mathbb{H}$. Moreover, let $\frac{d^-}{dy}$ and $\frac{d^+}{dy}$ denote the left and the right differential operators, respectively, that is,

$$\frac{d^-}{dy} h(y_0) := \lim_{\varepsilon \searrow 0} \frac{h(y_0) - h(y_0 - \varepsilon)}{-\varepsilon}$$

and

$$\frac{d^+}{dy} h(y_0) := \lim_{\varepsilon \searrow 0} \frac{h(y_0 + \varepsilon) - h(y_0)}{\varepsilon}.$$

PROPOSITION 3.1. *Assume that $y_1 \in (0,\infty)$ is such that $H_1(y_1) = W(y_1) < H_2(y_1)$. Also assume that the left and right derivatives $\frac{d^-}{dy} H_1$ and $\frac{d^+}{dy} H_1$ exist at $y_1$. Then*

(3.1) $$\frac{d^-}{dy} H_1(y_1) \geq \frac{d^-}{dy} W(y_1) \geq \frac{d^+}{dy} W(y_1) \geq \frac{d^+}{dy} H_1(y_1).$$

*Similarly, if $y_2 \in (0,\infty)$ is such that $H_2(y_2) = W(y_2)$ and $\frac{d^-}{dy} H_2$ and $\frac{d^+}{dy} H_2$ exist at $y_2$, then*

(3.2) $$\frac{d^-}{dy} H_1(y_2) \leq \frac{d^-}{dy} W(y_2) \leq \frac{d^+}{dy} W(y_2) \leq \frac{d^+}{dy} H_1(y_2).$$

PROOF. Since $W(y_1) = H_1(y_1)$, the first and the third inequality in (3.1) follow from $V \geq H_1$. Since $W(y_1) < H_2(y_1)$, we know that $W$ is concave in a neighborhood of $y_1$. From this, the second inequality follows.

The inequalities in (3.2) follow similarly. □

REMARK. Note that for the middle inequalities in (3.1) and (3.2) to hold, it is essential that $W(y_1) < H_2(y_1)$ and $H_1(y_2) < W(y_2)$, respectively. Indeed, (3.1) is, for example, not true at the point $y_1 = K_1$ if

$$H_1(y) = (y \wedge K_3 - K_2)^+$$

and

$$H_2(y) = (y - K_1)^+$$

for some constants $K_3 > K_2 > K_1 > 0$.

After a change of coordinates, Proposition 3.1 translates to the following smooth-fit principle. Note that, in line with the above results, no integrability conditions are assumed.

COROLLARY 3.2 (Smooth-fit principle). *Let $x_0 \in (0, \infty)$ and assume that $V(x_0) = g_i(x_0)$, where either $i = 1$ or $i = 2$. Assume also that $g_1(x_0) < g_2(x_0)$ and that $g_i$ is differentiable at $x_0$. Then also $V$ is differentiable at $x_0$ and*

$$\frac{d}{dx}V(x_0) = \frac{d}{dx}g_i(x_0).$$

**4. Existence of a saddle point.** According to Theorem 2.5, $\gamma^*$ is an optimal stopping time for the seller. It turns out, however, that

$$\tau^* := \inf\{t : V(X(t)) = g_1(X(t))\}$$

in general need not be optimal for the buyer; compare the examples in Section 5. A necessary condition for $(\tau^*, \gamma^*)$ to be a saddle point is that

$$\mathbb{P}(\tau^* < \infty) > 0,$$

or, equivalently, that the set

$$E_1 := \{x \in (0, \infty) : V(x) = g_1(x)\}$$

is nonempty. Indeed, $R_x(\infty, \infty) = 0$, and thus, $\tau^* = \infty$ cannot be optimal for the buyer (at least not if $g_1 \not\equiv 0$). Below we give an analytical criterion in terms of the differential operator

$$\mathcal{L} := \frac{\sigma^2}{2}\frac{\partial^2}{\partial x^2} + \mu\frac{\partial}{\partial x} - \beta,$$

ensuring that the set $E_1$ is empty. To this end, we restrict the class of payoff functions by requiring some additional regularity conditions.

HYPOTHESIS 4.1. Let $D = \{a_1, \ldots, a_n\}$, where $n \in \mathbb{N}$ and $a_i$ are positive real numbers with $a_1 < a_2 < \cdots < a_n$. Suppose that $g_1$ is a continuous function on $(0, \infty)$ such that $g_1'$ and $g_1''$ exist and are continuous on $(0, \infty) \setminus D$ and that the limits

$$g_1'(a_i\pm) := \lim_{x \to a_i\pm} g_1'(x), \qquad g_1''(a_i\pm) := \lim_{x \to a_i\pm} g_1''(x)$$

exist and are finite.

PROPOSITION 4.2. *Assume that the function $g_1$ satisfies Hypothesis 4.1 and that $g_2 > g_1$ on some open interval $\mathcal{I} \subset (0, \infty)$. If $\mathcal{L}g_1$ is a nonzero nonnegative measure on $\mathcal{I}$, then $V(x) > g_1(x)$ for every $x \in \mathcal{I}$. Thus, if $\mathcal{I} = (0, \infty)$, then the set $E_1$ is empty, and consequently, $\tau^*$ is not optimal for the buyer (provided $g_1 \not\equiv 0$).*

*Similarly, if $\mathcal{L}g_2$ is a nonzero nonpositive measure on $\mathcal{I}$, then $V(x) < g_2(x)$ for all $x \in \mathcal{I}$.*

REMARK. That $\mathcal{L}g_1$ is a nonnegative measure on $\mathcal{I}$ means that $\mathcal{L}g_1(x) \geq 0$ for all $x \in \mathcal{I} \setminus D$ and $g_1'(a-) \leq g_1'(a+)$ for all $a \in \mathcal{I} \cap D$. That $\mathcal{L}g_1$ is a nonzero nonnegative measure on $\mathcal{I}$ means that at least one of these inequalities is strict.

PROOF OF PROPOSITION 4.2. Fix $x \in \mathcal{I}$ and choose $l, r \in \mathcal{I}$ with $l < x < r$ so that $\mathcal{L}g_1$ is a nonzero nonnegative measure on $(l, r) \subset \mathcal{I}$. According to Theorem 2.5, $V(x) = \sup_\tau R_x(\tau, \gamma^*)$ and thus

(4.1) $\quad V(x) \geq R_x(\tau_{l,r}, \gamma^*) \geq \mathbb{E}_x(e^{-\beta(\tau_{l,r} \wedge \gamma^*)} g_1(X(\tau_{l,r} \wedge \gamma^*)))$.

Note that if $\mathbb{P}_x(\gamma^* < \tau_{l,r}) > 0$, then the second inequality in (4.1) is strict. Because $g_1$ satisfies Hypothesis 4.1, the Itô–Tanaka formula (see Theorem 3.7.1, page 218 in [12]) gives

$$\begin{aligned}
&\mathbb{E}_x(e^{-\beta(\tau_{l,r} \wedge \gamma^*)} g_1(X(\tau_{l,r} \wedge \gamma^*))) \\
&\quad = g_1(x) + \mathbb{E}_x\left(\int_0^{\tau_{l,r} \wedge \gamma^*} e^{-\beta s} \mathcal{L}g_1(X(s))\, ds\right) \\
&\qquad + \sum_{a_i \in (l,r)} (g_1'(a_i+) - g_1'(a_i-)) \mathbb{E}_x\left(\int_0^{\tau_{l,r} \wedge \gamma^*} e^{-\beta s}\, dL^i(s)\right),
\end{aligned}$$

where $L^i$ is the local time of $X$ at $a_i$. Now, since $\mathcal{L}g_1$ is nonnegative on $(l, r)$, we find that

$$\mathbb{E}_x(e^{-\beta(\tau_{l,r} \wedge \gamma^*)} g_1(X(\tau_{l,r} \wedge \gamma^*))) \geq g_1(x).$$

Moreover, if $\gamma^* \geq \tau_{l,r}$ a.s., then this inequality is strict. Indeed, since $\mathcal{L}g_1$ is a nonzero nonnegative measure on $(l, r)$, we have that either $\mathcal{L}g_1(y) > 0$ for some $y \in (l,r)$, where $g_1$ is differentiable (implying that the middle term is strictly positive), or $g_1'(a_i+) > g_1'(a_i-)$ for some $a_i \in (l,r)$ (implying that the last term is strictly positive). Thus, in view of (4.1), we have $V(x) > g_1(x)$, which finishes the proof of the first part of the proposition.

As for the second claim, by Proposition 4.4 in [6] [note that it is also valid for contracts, functions of the type (2.7)], we may replace (4.1) with

$$V(x) \leq \sup_{\tau} R_x(\tau, \gamma_{l,r}) = R_x(\hat{\tau}, \gamma_{l,r})$$

for some stopping time $\hat{\tau}$. The proof now follows as above. □

Below we provide conditions under which $\tau^*$ is optimal for the buyer. Following [1] and [6], the conditions are expressed in terms of the two quantities

$$l_0 := \limsup_{x \to 0} \frac{g_1(x)}{\varphi(x)} \quad \text{and} \quad l_\infty := \limsup_{x \to \infty} \frac{g_1(x)}{\psi(x)}.$$

PROPOSITION 4.3. *Assume that both $l_0$ and $l_\infty$ are finite. Also assume that the nonnegative local martingales $e^{-\beta t}\varphi(X(t))$ and $e^{-\beta t}\psi(X(t))$ satisfy*

$$(4.2) \quad \mathbb{E}_x\left(\sup_{0 \leq s \leq t} e^{-\beta s}\varphi(X(s))\right) < \infty \quad \text{and} \quad \mathbb{E}_x\left(\sup_{0 \leq s \leq t} e^{-\beta s}\psi(X(s))\right) < \infty$$

*for all times $t$. Then the process $e^{-\beta t \wedge \tau^*}V(X(t \wedge \tau^*))$ is a sub-martingale.*

PROOF. We know from Theorem 2.5 that $V/\varphi$ is $F$-convex in all intervals where $V > g_1$. Arguing as in the proof of Proposition 5.1 in [6], it can therefore be shown that $Z(t) := e^{-\beta t \wedge \tau^*}V(X(t \wedge \tau^*))$ is a sub-martingale, provided

$$E_x\left(\sup_{0 \leq s \leq t} Z(s)\right) < \infty$$

(this is needed for the use of Fatou's lemma). From the results in [6] (compare Propositions 5.4 and 5.12 of that paper) we know that

$$\limsup_{x \to 0} \frac{V(x)}{\varphi(x)} = l_0 \quad \text{and} \quad \limsup_{x \to \infty} \frac{V(x)}{\psi(x)} = l_\infty.$$

Thus, there exist constants $C$ and $D$ with

$$V(x) \leq C\varphi(x) + D\psi(x)$$

for all $x \in (0, \infty)$. From the assumption (4.2), it therefore follows that $\sup_{0 \leq s \leq t} Z(s)$ is integrable, which completes the proof. □

REMARK. Without the assumption (4.2), Proposition 4.3 would not be true. Also note that to show that the process $e^{-\beta t \wedge \gamma^*} V(X(t \wedge \gamma^*))$ is a super-martingale, neither the finiteness of $l_0$ and $l_\infty$ nor the condition (4.2) is needed.

The following two results may be viewed as the game versions of Proposition 5.13 and 5.14 in [6].

THEOREM 4.4. *Assume (4.2) and that*

(4.3) $$l_0 = l_\infty = 0.$$

*Then $(\tau^*, \gamma^*)$ is a saddle point.*

PROOF. From Proposition 4.3, it follows that

$$\begin{aligned} V(x) &\leq \mathbb{E}_x(e^{-\beta(t \wedge \tau^* \wedge \gamma)} V(X(t \wedge \tau^* \wedge \gamma))) \\ &\leq \mathbb{E}_x(e^{-\beta \tau^*} V(X(\tau^*)) \mathbb{1}_{\{\tau^* \leq t \wedge \gamma\}} + e^{-\beta \gamma} V(X(\gamma)) \mathbb{1}_{\{\gamma < t \wedge \tau^*\}}) \\ &\quad + \mathbb{E}_x(e^{-\beta t} V(X(t)) \mathbb{1}_{\{t \leq \tau^* \wedge \gamma\}}) \end{aligned}$$

for any stopping time $\gamma$. We first prove that the last term converges to zero when $t$ tends to $+\infty$. To do this, recall that the assumption (4.3) implies

$$\lim_{x \to 0} \frac{V(x)}{\varphi(x)} = \lim_{x \to \infty} \frac{V(x)}{\psi(x)} = 0.$$

Thus, given a constant $\delta > 0$, there exists a constant $M$ such that $V(x) \leq \delta \varphi(x) + \delta \psi(x) + M$ for all $x$. Using the fact that $e^{-\beta t} \varphi(X(t))$ and $e^{-\beta t} \psi(X(t))$ are nonnegative local martingales, and hence supermartingales, we find

$$\begin{aligned} \mathbb{E}_x(e^{-\beta t} V(X(t)) \mathbb{1}_{\{t \leq \tau^* \wedge \gamma\}}) &\leq M e^{-\beta t} + \delta \mathbb{E}_x e^{-\beta t} \varphi(X(t)) + \delta \mathbb{E}_x e^{-\beta t} \psi(X(t)) \\ &\leq M e^{-\beta t} + \delta \varphi(x) + \delta \psi(x). \end{aligned}$$

Since $\delta$ can be chosen arbitrarily, we conclude the first step. Next, the monotone convergence theorem yields

$$\begin{aligned} V(x) &\leq \lim_{t \to \infty} \mathbb{E}_x(e^{-\beta \tau^*} V(X(\tau^*)) \mathbb{1}_{\{\tau^* \leq t \wedge \gamma\}} + e^{-\beta \gamma} V(X(\gamma)) \mathbb{1}_{\{\gamma < t \wedge \tau^*\}}) \\ &\leq \lim_{t \to \infty} \mathbb{E}_x(e^{-\beta \tau^*} g_1(X(\tau^*)) \mathbb{1}_{\{\tau^* \leq t \wedge \gamma\}} + e^{-\beta \gamma} g_2(X(\gamma)) \mathbb{1}_{\{\gamma < t \wedge \tau^*\}}) \\ &= \mathbb{E}_x(e^{-\beta \tau^*} g_1(X(\tau^*)) \mathbb{1}_{\{\tau^* \leq \gamma\}} + e^{-\beta \gamma} g_2(X(\gamma)) \mathbb{1}_{\{\gamma < \tau^*\}}) \\ &= R_x(\tau^*, \gamma), \end{aligned}$$

that is, $\tau^*$ is optimal for the buyer. This completes the proof. □

THEOREM 4.5. *Assume* (4.2) *and that $l_0$ and $l_\infty$ are both finite. Then, the pair $(\tau^*, \gamma^*)$ is a saddle point for arbitrary starting point if and only if*

$$\left\{ \begin{array}{c} \textit{there is no } l>0 \textit{ such that} \\ g_1(x) < V(x) \textit{ for all } x \leq l \\ \textit{if } l_0 > 0 \end{array} \right\} \quad \textit{and} \quad \left\{ \begin{array}{c} \textit{there is no } r>0 \textit{ such that} \\ g_1(x) < V(x) \textit{ for all } x \geq r \\ \textit{if } l_\infty > 0 \end{array} \right\}.$$

PROOF. If $l_0 = l_\infty = 0$, then the result follows from Theorem 4.4. Therefore, we assume that $l_\infty > 0$ (the case $l_0 > 0$ can be treated similarly).

To prove the sufficiency of the condition, fix a starting point $x \in (0, \infty)$. If $V(x) = g_1(x)$, then $\tau^* = 0$ is clearly optimal for the buyer, and thus, we are finished. If $V(x) > g_1(x)$, let $I := (a, b) \subset (0, \infty)$ be a maximal interval containing $x$ such that $V > g_1$ in $I$. Note that

$$\tau^* = \inf\{t : X(t) \notin I\},$$

and that $b < \infty$ by assumption. Moreover, given $\delta > 0$, there exists a constant $M$ such that $V \leq M + \delta\varphi$ in $I$. Indeed, if $a > 0$, then $V$ is bounded in $I$, and if $a = 0$, then $l_0 = 0$ by assumption. Thus, proceeding analogously as in the proof of Theorem 4.4, we obtain

$$\begin{aligned} V(x) &\leq \lim_{t \to \infty} \mathbb{E}_x(e^{-\beta(t \wedge \tau^* \wedge \gamma)} V(X(t \wedge \tau^* \wedge \gamma))) \\ &\leq \lim_{t \to \infty} \mathbb{E}_x(e^{-\beta\tau^*} V(X(\tau^*)) \mathbb{1}_{\{\tau^* \leq t \wedge \gamma\}} + e^{-\beta\gamma} V(X(\gamma)) \mathbb{1}_{\{\gamma < t \wedge \tau^*\}}) \\ &\quad + \delta\varphi(x) + \lim_{t \to \infty} M e^{-\beta t} \\ &\leq R_x(\tau^*, \gamma) + \delta\varphi(x) \end{aligned}$$

for a stopping time $\gamma$. Since $\delta$ is arbitrary, this shows that $\tau^*$ is optimal for the buyer.

Conversely, assume that $(\tau^*, \gamma^*)$ is a saddle point for each starting point $x$ and that $V(x) > g_1(x)$ for $x \geq r$. Then, for $x \geq r$, the stopping time $\tau^* \geq \tau_r$ a.s. The definition of a saddle point and the optional sampling theorem applied to the nonnegative supermartingale $e^{-\beta t} V_\infty(X_t)$, where $V_\infty$ is the perpetual American option value as defined in (2.14), give

$$\begin{aligned} V(x) &= R_x(\tau^*, \gamma^*) \\ &\leq R_x(\tau^*, \infty) = \mathbb{E}_x(e^{-\beta\tau^*} g_1(X(\tau^*))) \\ &\leq \mathbb{E}_x(e^{-\beta\tau^*} V_\infty(X(\tau^*))) \\ &\leq \mathbb{E}_x(e^{-\beta\tau_r} V_\infty(X(\tau_r))) \\ &= \frac{\varphi(x)}{\varphi(r)} V_\infty(r), \end{aligned}$$

where we in the last equation used (2.6) in [6]. Proposition 5.4 in [6] then implies that

$$l_\infty = \limsup_{x\to\infty} \frac{V(x)}{\psi(x)} \leq \frac{V_\infty(r)}{\varphi(r)} \lim_{x\to\infty} \frac{\varphi(x)}{\psi(x)} = 0,$$

which contradicts $l_\infty > 0$. □

**5. Two examples of game options.** In this section we study two examples motivated by applications in finance. In both examples we assume that $\mu(x) = \beta x$, where $\beta$ is the discounting rate. Thus, the diffusion $X$ solves

$$dX(t) = \beta X(t)\, dt + \sigma(X(t))\, dW(t),$$

and $V$ may be interpreted as the arbitrage free price of a game option written on a nondividend paying stock; compare [14]. Note that the functions $\psi$ and $\varphi$ are given (up to multiplication with a positive constant) by

$$\psi(x) = x$$

and

(5.1) $$\varphi(x) = x \int_x^\infty \frac{1}{u^2} \exp\left\{-\int_1^u \frac{2\beta z}{\sigma^2(z)}\, dz\right\} du.$$

5.1. *The game version of a call option.* In this subsection we study the game version of a call option, that is,

$$g_1(x) = (x - K)^+ \quad \text{and} \quad g_2(x) = (x - K)^+ + \varepsilon$$

for some positive constants $K$ and $\varepsilon$. If $\varepsilon \geq K$, then one can show that the game option reduces to an ordinary perpetual American call option. Therefore, we consider the case with $\varepsilon < K$.

The functions $H_i := (\frac{g_i}{\varphi}) \circ F^{-1}$, $i = 1, 2$, are given by

$$H_1(y) = \left(y - \frac{K}{\varphi(F^{-1}(y))}\right)^+$$

and

$$H_2(y) = \left(y - \frac{K}{\varphi(F^{-1}(y))}\right)^+ + \frac{\varepsilon}{\varphi(F^{-1}(y))}.$$

First we claim that the function

$$w(y) := \frac{1}{\varphi(F^{-1}(y))}$$

is concave. To see this, note that by letting $y = F(x)$, we find that

$$w(y) = \frac{1}{\varphi(F^{-1}(y))} = \frac{1}{\varphi(x)} = \frac{F(x)}{x} = \frac{y}{F^{-1}(y)},$$

where we have used $F(x) = x/\varphi(x)$. Straightforward calculations yield that

$$w''(y) - \frac{\varphi''(x)}{\varphi^3(x)(F'(x))^2}.$$

Using (5.1), one can check that $\varphi''(x) \geq 0$, so it follows that $w$ is concave. Since $w$ is concave, $H_1$ is 0 on $(0, F(K))$ and convex in $(F(K), \infty)$, and $H_2$ is concave in $(0, F(K))$ and convex in $(F(K), \infty)$. This, together with the easily checked facts

$$\lim_{y \to \infty} \frac{H_1(y)}{y} = 1, \qquad H_2'(y) < 1$$

and

$$H_2'(F(K)+) = \frac{\varepsilon}{K} + \frac{(K-\varepsilon)F(K)}{K^2 F'(K)} > \frac{\varepsilon}{K} = \frac{H_2(F(K))}{F(K)},$$

implies that the smallest function $W$ in $\mathbb{H}$ is given by

$$W(y) = \begin{cases} \dfrac{\varepsilon y}{K}, & \text{if } y \in (0, F(K)], \\ H_2(y), & \text{if } y \in (F(K), \infty). \end{cases}$$

In the usual coordinates this means that the value $V$ of the game version of a call option written on a no-dividend paying stock is

$$V(x) = \begin{cases} \dfrac{\varepsilon x}{K}, & \text{if } x \in (0, K], \\ x - K + \varepsilon, & \text{if } x \in (K, \infty). \end{cases}$$

According to Theorem 2.5, an optimal stopping time for the seller is given by

$$\gamma^* := \inf\{t : X(t) \geq K\}.$$

Also note that the corresponding stopping time $\tau^* = \infty$ is not optimal for the buyer.

5.2. *An example in which convexity is lost.* In this subsection we consider another possible generalization of the American call option. More precisely, let

$$g_1(x) = (x - K)^+ \quad \text{and} \quad g_2(x) = C(x - K)^+$$

for some constant $C > 1$. Moreover, assume for simplicity that the diffusion $X$ is a geometric Brownian motion, that is, that

$$dX(t) = \beta X(t)\,dt + \sigma X(t)\,dW(t)$$

for some constant $\sigma > 0$. Then the functions $\psi$ and $\varphi$ are given by

$$\psi(x) = x \quad \text{and} \quad \varphi(x) = x^{-2\beta/\sigma^2},$$

and the functions $H_i$, $i = 1, 2$, are given by

$$H_1(y) = (y - Ky^{2\beta/(2\beta+\sigma^2)})^+ \quad \text{and} \quad H_2(y) = C(y - Ky^{2\beta/(2\beta+\sigma^2)})^+.$$

We need to consider two different cases.

5.2.1. *Case* 1. First assume that $C \geq 1 + 2\beta/\sigma^2$. Then it is straightforward to check that $W(y) = (y - K^{(2\beta+\sigma^2)/\sigma^2})^+$, that is, the value $V$ of the option is given by

$$V(x) = \varphi(x)W(F(x)) = (x - K^{(2\beta+\sigma^2)/\sigma^2} x^{-2\beta/\sigma^2})^+.$$

Moreover, Theorem 2.5 tells us that $\gamma^* := \inf\{t : X(t) \leq K\}$ is an optimal stopping time for the seller.

5.2.2. *Case* 2. Now assume that $1 < C < 1 + 2\beta/\sigma^2$. Then one can check that

$$W(y) = \begin{cases} H_2(y), & \text{if } y \in (0, y'), \\ H_2(y') + y - y', & \text{if } y \in [y', \infty), \end{cases}$$

where $y'$ is given by

$$y' = \left(\frac{2\beta CK}{(2\beta + \sigma^2)(C - 1)}\right)^{(2\beta+\sigma^2)/\sigma^2}.$$

It follows that

$$V(x) = \begin{cases} C(x - K)^+, & \text{if } x \in (0, x'), \\ x - \dfrac{CK\sigma^2}{2\beta + \sigma^2}\left(\dfrac{x'}{x}\right)^{2\beta/\sigma^2}, & \text{if } x \in [x', \infty), \end{cases}$$

where

$$x' = \frac{2\beta CK}{(2\beta + \sigma^2)(C - 1)}.$$

According to Theorem 2.5, $\gamma^* := \inf\{t : X(t) \leq x'\}$ is optimal for the seller. As in the previous example, however, $\tau^* = \inf\{t : X(t) \leq K\}$ is not optimal for the buyer.

REMARK. The above example shows, perhaps surprisingly, that game options are not convexity preserving. More precisely, although both contract functions $g_1$ and $g_2$ are convex, the value of the game option need not necessarily be convex. This is in contrast to options of European and American style, both of which are known to be convexity preserving; compare, for example, [4] or [9] and the references therein.

Remark. The method to determine the value of an optimal stopping game used in this section is also used in [8]. In that paper the construction of the value using concave functions is shown to be valid under the assumption of the existence of a value and a saddle point of the form $(\tau^*, \gamma^*)$. In the present paper we start with the construction of a natural candidate for the value function (without knowing a priori that such a value function exists), and then we show that this function indeed has to be the value of the game. This allows us to weaken the assumptions under which a game is known to have a value. Also note that the integrability condition (1.4) is satisfied in neither of the two examples provided in this section.

**Acknowledgment.** The authors thank an anonymous referee for very helpful comments.

School of Mathematics  
University of Manchester  
Sackville Street  
Manchester M60 1QD  
United Kingdom  
E-mail: ekstrom@maths.manchester.ac.uk

GREMAQ UMR CNRS 5604  
Université des Sciences Sociales  
21 allée de Brienne  
31000 Toulouse  
France  
E-mail: stephane.villeneuve@univ-tlse1.fr